\documentclass[%
 aip,
cp,  
 amsmath,amssymb,
 reprint,%
]{revtex4-2}

\usepackage{graphicx}
\usepackage{dcolumn}
\usepackage{bm}

\usepackage[utf8]{inputenc}
\usepackage[T1]{fontenc}
\usepackage{mathptmx} 

\begin{document}

\title{A Computer Program for the Numerical Analysis of Economic Cycles
	Within the Framework of the Dubovsky Generalized Model}

\author{Makarov Danil} 
\email[Corresponding author: ]{danil.makarov.pk@yandex.ru}
\affiliation{Vitus Bering Kamchatka State University, Petropavlovsk-Kamchatskiy, Russia
}

\author{Parovik Roman}%
\email{romanparovik@gmail.com}
\affiliation{Vitus Bering Kamchatka State University, Petropavlovsk-Kamchatskiy, Russia
}

\date{\today} 

\begin{abstract}
The article proposes a computer program for calculating
economic crises according to the generalized mathematical model of S.V. Dubovsky.
This model is represented by a system of ordinary nonlinear differential
equations with fractional derivatives in the sense of Gerasimov-Caputo with
initial conditions. Furthermore, according to a numerical algorithm based on an
explicit nonlocal finite-difference scheme, oscillograms and phase trajectories
were constructed. It is shown that changing the orders of fractional derivatives
in the model can give rise to various modes, for example, damped modes with a
steady-state amplitude. It is concluded that the orders of fractional derivatives
are responsible for the intensity of the process.
\end{abstract}

\maketitle

\section{Introduction}

In the modern world, increasingly people need mathematical modeling of economic processes. This is due to the fact that when doing business, the mathematical
description of economic processes gives a quantitative and qualitative
presentation, which helps in further forecasting and making the right management decisions. One of the most important the the economic processes is economic crisis. Economic crises determine the economic well-being of citizens and the degree of social tension in the country. Back in the 1920s, the Soviet economist N.D. Kondratyev singled out long-term periodic fluctuations (waves) with a duration of 50-55 years in the economic time series \cite{Kond_1928}. 

The most complete mathematical description of modeling Kondrat'ev's cycles was carried out in the works of S.V. Dubovsky \cite{Dubovsky_1995}. In this paper, a mathematical model was proposed that generalizes the well-known Dubovsky model in the case of taking into account the effects of memory in the economic system and is a logical continuation of this
work \cite{Makarov_2016}. Memory effects are described using the theory of fractional calculus, namely, fractional derivatives \cite{Nakh_2003}.

\section{Mathematical model}
The generalized model of Kondratyev's cycles can be represented as a system of equations.

\begin{equation}
\left\{ \begin{array}{l}
\partial _{0t}^\alpha x\left( \tau  \right) =  - \lambda nx\left( t
\right)\left( {x\left( t \right) - 1} \right)\left( {y\left( t \right) - {y^*}}
\right)\\
\partial _{0t}^\beta x\left( \tau  \right) = n\left( {1 - n} \right){y^2}\left(
t \right)\left( {x\left( t \right) - {x^*}} \right) + f\left( t \right)\\
x\left( 0 \right) = a,y(0) = b.
\end{array} \right.
\label{eq1}
\end{equation}
where $\partial _{0t}^\alpha x\left( \tau  \right) = \frac{1}{{\Gamma \left( {1
			- \alpha } \right)}}\int\limits_0^t {\frac{{\dot x\left( \tau  \right)d\tau
	}}{{{{\left( {t - \tau } \right)}^\alpha }}}} ,\partial _{0t}^\beta x\left( \tau 
\right) = \frac{1}{{\Gamma \left( {1 - \beta } \right)}}\int\limits_0^t
{\frac{{\dot x\left( \tau  \right)d\tau }}{{{{\left( {t - \tau } \right)}^\beta
}}}} $, $0 < \alpha ,\beta  < 1$ --- fractional derivatives in the sense of
Gerasimov-Caputo; $\Gamma \left( x \right)$ --- Euler's gamma function; $x\left( t
\right)$ --- the effectiveness of new technologies; $y\left( t \right)$  ---
efficiency of return on assets; $x^*$and $y^*$ --- equilibrium stationary solution
of system (1); $n$  ---  the rate of accumulation; $\lambda $  --- coefficient,
which is determined from the statistics of the time series; $f\left( t \right)$ --- external impact on the economic
system; $t \in \left[ {0,T} \right]$ --- time coordinate, $T$ --- process
simulation time; $a$  and $b$  --- initial conditions, given constants.

Note that the nonlinear system (1) in the case of the values of the parameters 
$\alpha  = \beta  = 1$ and $f\left( t \right) = 0$ goes over to the model of S.V.
Dubovsky \cite{Dubovsky_1995}. The solution to the nonlinear system (1) will be sought using
numerical methods - finite difference schemes. Let us split the time interval
$\left[ {0,T} \right]$ into $N$ equal parts with a step $\tau  = {T
	\mathord{\left/
		{\vphantom {T N}} \right.
		\kern-\nulldelimiterspace} N}$.
	
The approximation of the fractional derivatives in the equation is carried out
according to \cite{Makarov_2016}. Then the system will be written in the finite-difference
formulation in the form of the system.	

\begin{equation}
\left\{ \begin{array}{l}
{x_0} = a,{y_0} = b,\\
{x_1} = {x_0}\left( {1 - \frac{{\lambda n}}{A}\left( {{x_0} - 1} \right)\left(
	{{y_0} - {y^ * }} \right)} \right),{y_1} = {y_0}\left( {1 + \frac{{n\left( {1 -
				n} \right)}}{B}{y_0}\left( {{x_0} - {x^ * }} \right)} \right),j = 0,\\
{x_{j + 1}} = {x_j}\left( {1 - \frac{{\lambda n}}{A}\left( {{x_j} - 1}
	\right)\left( {{y_j} - {y^ * }} \right)} \right) - \sum\limits_{k = 1}^{j - 1}
{{p_k}\left( {{x_{j - k + 1}} - {x_{j - k}}} \right),} \\
{y_{k + 1}} = {y_j}\left( {1 + \frac{{n\left( {1 - n} \right)}}{B}{y_j}\left(
	{{x_j} - {x^ * }} \right)} \right) - \sum\limits_{k = 1}^{j - 1} {{q_k}\left(
	{{y_{j - k + 1}} - {y_{j - k}}} \right) + {f_j},j = 1,...,N - 1,}
\end{array} \right.
\label{eq2}
\end{equation}
where $A = \frac{{{\tau ^{ - \alpha }}}}{{\Gamma \left( {2 - \alpha }
		\right)}},B = \frac{{{\tau ^{ - \beta }}}}{{\Gamma \left( {2 - \beta }
		\right)}},{p_k} = {\left( {1 + k} \right)^{1 - \alpha }} - {k^{1 - \alpha
}},{q_k} = {\left( {1 + k} \right)^{1 - \beta }} - {k^{1 - \beta }}$.

Let us investigate solution (2) depending on various values of the fractional
parameters \textit{$\alpha{}$} and \textit{$\beta{}$}, and construct the phase
trajectories. In this paper, we do not dwell on the questions of stability or
convergence of the explicit finite difference scheme (2).

\section{Environment and programming language}

The high-level Python programming language and PyCharm development environment
were chosen for the development of the application. Python supports structured,
object-oriented, functional, imperative, and aspect-oriented programming. The
main architectural features are dynamic typing, automatic memory management, full
introspection, an exception handling mechanism, support for multithreaded
calculations, and high-level data structures. The division of programs into
modules is supported, which, in turn, can be combined into packages.

PyCharm is an integrated development environment for the Python programming
language. Provides code analysis tools, a graphical debugger, a unit test runner,
and supports Django web development.

PyCharm features \cite{Lutz_2014}:

1. Static code analysis, syntax highlighting and errors.

2. Navigation through the project and source code.

3. Display of the project file structure, quick transition between files,
classes, methods, and method uses.

4. Refactoring: renaming, retrieving a method, introducing a variable,
introducing a constant, raising and lowering.

5. Tools for web development using the Django framework.

6. Built-in debugger for Python.

7. Built-in unit testing tools.

8. Development of GoogleAppEngine.

9. Support for version control systems: common user interface for Mercurial, Git, Subversion, Perforce, and CVS with support for changelog and merge.

The Kivy library was used to create a beautiful graphical interface.

\section{Simulation results}

We take the modeling parameters from [3]: \textit{x}$^{*}$= 1.3,
\textit{y}$^{*}$= 0.5, \textit{n} = 0.2, \textit{$\lambda{}$} = 2.25,
\textit{x}(0) = 1.35, \textit{y}(0) = 0.5, \textit{T}=250, \textit{$\alpha{}$} =
\textit{$\beta{}$} = 1 (Figure 1). In this case, you can see that the phase
trajectory has an ellipsoidal closed shape, the equilibrium state of the state of
the system is called the center. The vibration amplitude is constant

\begin{figure}[h!]
\centering
\includegraphics[scale=1.1]{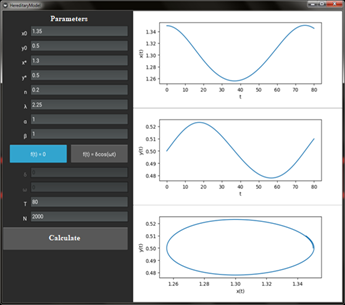}
\caption{Calculated curves and phase trajectory in the case $f\left(t\right) = 0$}
\end{figure}

Let us add to the system the influence of the external periodic influence $f\left(t\right) = \delta\cos\left(\omega t\right)$ --- investment cycles. The value of the parameters is $\delta = 0.01$ and $\omega = 1$ (Figure 2). It can be concluded that the external periodic impact gives a cycle with a period of about 7 years, and the main cycle is 60 years, which corresponds to the upper limit of the Kondratyev cycle \cite{Dubovsky_1995}. This combined model is the most flexible in describing economic crises.

\begin{figure}[h!]
\centering
\includegraphics[scale=1.2]{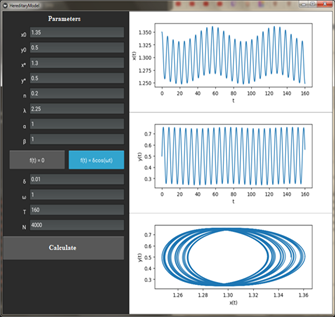}
\caption{Calculated curves and phase trajectory in the case $f\left(t\right) = \delta\cos\left(\omega t\right)$}
\end{figure}
	
Consider the case when $f\left(t\right)=0, \alpha=0.8$ and $\beta=1$, and the other parameters remain unchanged (Figure 3). The calculated curve shows that the oscillation process is damped and the phase trajectory is not closed. The equilibrium position is called a stable focus.

\begin{figure}[h!]
\centering
\includegraphics[scale=1]{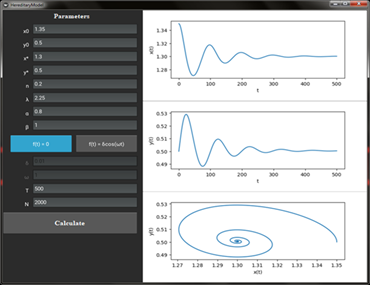}
\caption{Calculated curves and phase trajectory in the case of $\alpha=0.8, \beta=1$ and $f\left(t\right)=0$}
\end{figure}

Add the external influence $f\left(t\right)=\delta\cos\left(\omega t\right), \alpha=0.8, \beta=0.6, \delta=0.5, \omega=2$. (Figure 4). First, the amplitude of oscillations increases, and then it enters a constant mode, this can be seen on the phase trajectory, which eventually reaches a constant mode or limit cycle, which can be used in the study of Kondratiev cycles.

\begin{figure}[h!]
\centering
\includegraphics[scale=1]{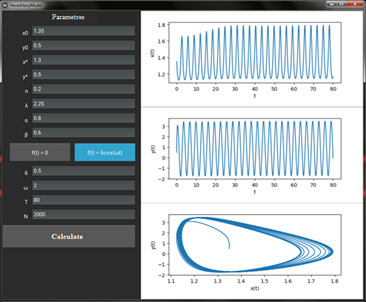}
\caption{Calculated curve and phase trajectory in the case of $f\left(t\right)=\delta\cos\left(\omega t\right), \alpha=0.8, \beta=0.6, \delta=0.5, \omega=2$}
\end{figure}

Consider 2 more cases $\alpha=0.8, \beta=0.8$ and $f\left(t\right) = \delta\cos\left(\omega t\right), \delta=0.5, \omega=2$ (Figure 5) and $\alpha=0.1, \beta=0.1$ and $f\left(t\right) = \delta\cos\left(\omega t\right), \delta=0.5, \omega=2$ (Figure 6). In these cases, the phase trajectories reach the limit cycle.

\begin{figure}[h!]
\centering
\includegraphics[scale=1]{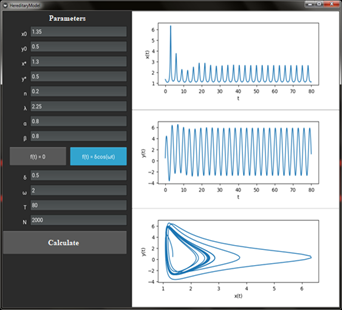}
\caption{ Calculated curve and phase trajectory in the case of $\alpha=0.8, \beta=0.8$ and $f\left(t\right) = \delta\cos\left(\omega t\right), \delta=0.5, \omega=2$}
\end{figure}

\begin{figure}[h!]
\centering
\includegraphics[scale=1]{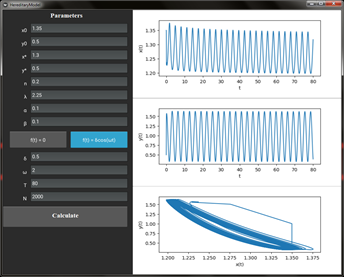}
\caption{Calculated curve and phase trajectory in the case of $\alpha=0.1, \beta=0.1$ and $f\left(t\right) = \delta\cos\left(\omega t\right), \delta=0.5, \omega=2$}
\end{figure}
	
\section{Conclusion}

During the work, a hereditary dynamic system was modeled that simulates the
economic cycles of Dubovsky, which takes into account the effects of memory in
the economic system. Some cases with and without external influence on the system are presented. Moreover, if we introduce fractional derivatives into the system, this will lead to damping processes, but if there is an external periodic influence in the system, then the system enters the limit cycle, which can be considered one or another economic cycle.

According to the research results, we see that the double orders of derivatives in system (1) are responsible for the intensity of the economic process and, as shown in the article \cite{Parovik_2019}, are related to its quality factor.

\begin{acknowledgments}
The work was performed within the framework of the research project of Vitus Bering KamSU "Mathematical model of Kondratiev's long waves taking into account heredity" No. AAAA-A20-120021190003-1.
\end{acknowledgments}


\end{document}